# A NOTE ON MAXIMUM AND MINIMUM STABILITY OF CERTAIN DISTRIBUTIONS

S. SATHEESH

and

N. UNNIKRISHNAN NAIR

*Cochin University of Science and Technology*
*Cochin - 682 022, India.*
ssatheesh@sancharnet.in

***ABSTRACT*:** In the context of stability of extremes of a random variable *X* with respect to a discrete random variable *N*, we discuss relations between *X* and *N* and extend the notions when *X* is discrete.

## 1. INTRODUCTION

Let $X>0$ be a continuous random variable (r.v) with distribution function (d.f) $F$, and $N>0$ be an integer valued r.v independent of $X$ with probability generating function (PGF) $Q(s)$. By Voorn (1987) $F$ is maximum stable with respect to (w.r.t) $N$ ($F$ is $N$-max stable) if,

$$Q[F(x)] = F(cx), \text{ for all } x>0 \text{ and some } c>0. \tag{1.1}$$

Similarly $F$ is minimum stable w.r.t $N$ ($F$ is $N$-min stable) if,

$$Q[R(cx)] = R(x), \text{ for all } x>0 \text{ and some } c>0, \tag{1.2}$$

where $R(x) = 1-F(x)$. This note addresses four problems. (i) when $X$ is exponential (ii) non-geometric laws for $N$ (iii) finding $N$ for the stability of a given $F$ and (iv) extending the notion to a discrete r.v $X$. Proofs of some of the results only are given. We will use $I_1$ for $\{1,2,\ldots\}$ and $I_0$ for $\{0,1,2,\ldots\}$. Notice that in (1.1) or (1.2) $c \in (0,1)$.

## 2. CHARACTERIZATIONS

**Proposition.2.1** Exponential law is N-max stable iff $N$ is Sibuya($v$) with PGF

$$Q(s) = 1-(1-s)^v, \ 0<v<1,$$

and $c = v$. It is N-min stable iff $N$ is degenerate at $k >1$ integer, $c = 1/k$.

**Proposition.2.2** $F$ is Sibuya($v$)-max stable iff it is semi-Weibull($p,\alpha$) with

$$F(x) = 1 - e^{-\psi(x)}, \ x>0, \ p = v = c^\alpha, \text{ and } \psi(x) \text{ satisfies}$$

$$\psi(x) = \psi(p^{1/\alpha} x)/p \quad \text{for all } x >0, \text{ some } 0<p<1, \text{ and } \alpha >0. \tag{2.1}$$

**Proof.** Setting $\psi(x) = -\log R(x)$, $Q[F(x)] = 1 - e^{-v\psi(x)}$.

Sibuya-max stability of $F$ implies, $1 - e^{-v\psi(x)} = 1 - e^{-\psi(cx)}$.

Hence $F$ is semi-Weibull($p,\alpha$), $p = v = c^\alpha$. Converse is easy.

**Proposition.2.3** $F$ is min stable w.r.t Harris($a,k$) law with PGF

$$Q(s) = s/[a-(a-1)s^k]^{1/k}, \ k \in I_1, \ a >1$$

iff it is generalized semi-Pareto($p,\alpha,\beta$), $p = 1/a = c^\alpha$, $\beta =1/k$ with

$$F(x) = 1-[1+ \psi(x)]^{-\beta}, \ x>0, \ \beta >0, \ \psi(x) \text{ satisfies (2.1)}.$$

**Proof.** With $\psi(x) = [1-R^k(s)]/R^k(s)$, $k \in I_1$, $F(x) = 1- [1+ \psi(x)]^{-1/k}$.

Harris-min stability of $F$ implies, $[1+ a\psi(cx)]^{-1/k} = [1+ \psi(x)]^{-1/k}$.

Hence $F$ is generalized semi-Pareto($1/a,\alpha,1/k$). Converse is easy.

This includes all Pareto types and generalizes the result on semi-Pareto laws by Pillai and Sandhya (1996) as Harris($a$,1) is geometric($1/a$) on $I_1$. Also, $F$ is min stable w.r.t a degenerate law at $k >1$ integer iff it is semi-Weibull($p,\alpha$), $pk = 1$.

Under (1.1): $Q[F(x)] = F(cx) = F_c(x)$. When $x = F^{-1}(s)$, $0<s<1$, $F_c[F^{-1}(s)] = Q(s)$. Similarly, $R[R_c^{-1}(s)] = Q(s)$ under (1.2). These emphasize the necessity of the $Q$ for the N-max/min stability of a given $F$. Sreehari (1995) has the relation between $F$ and the Laplace transform (LT) of $N$ for max stability. A study of N-sum stability motivated us. Here the relation between the parameters of $Q$ and $F$ and their ranges can be found by:

**Lemma.2.1** If $Q(s)$ is a PGF, then $Q(s^t)$ is a PGF iff $t > 0$ is an integer.

**Proof.** Since, $Q(s^t)$ is a power series in $s$ (Feller, 1971, 223) iff $t \in I_1$.

**Proposition.2.4** Generalized semi-Pareto$(p,\alpha,\beta)$ is N-min stable iff $N$ is Harris$(a,k)$, $a = 1/p$, $\beta = 1/k$, $c = p^{1/\alpha}$.

**Proof.** Here $R(x) = [1+\psi(p^{1/\alpha} x)/p]^{-\beta}$ and $R_c^{-1}(s) = \psi^{-1}(u)/c$, $u = (1-s^{1/\beta})/s^{1/\beta}$.

When $c = p^{1/\alpha}$, $R[R_c^{-1}(s)] = s/[a - (a-1)s^{1/\beta}]^\beta$, $a = 1/p$.

This is a PGF iff $1/\beta = k \in I_1$ implying $N$ is Harris$(a,k)$. Converse is easy.

**Proposition.2.5** Semi-Weibull$(p,\alpha)$ is N-max stable iff $N$ is Sibuya$(p)$ and N-min stable iff $N$ is degenerate at $1/p > 1$ integer. In both cases $c^\alpha = p$.

**Proposition.2.6** $F(x) = [1 + x^{-\alpha}]^{-1/k}$, $x > 0$, $k \in I_1$, $\alpha > 0$ (extended log-logistic law of Voorn (1987)) is N-max stable iff $N$ is Harris$(c^{-\alpha},k)$.

Finally, when $X$ is discrete: Let $\{m(j)\}$ is the realizations of a LT at $j \in I_0$. Now,

$$P\{X<j\} = F(j) = 1 - m(j), \quad j \in I_0 \tag{2.2}$$

is the d.f of a mixture of geometric$(p)$ laws on $I_0$ by Satheesh and Sandhya (1997). As $m(cs)$, $c>0$, is again a LT the function $G(j) = 1 - m(cj)$, is a d.f. Since $F(cj) = G(j)$, definitions of N-max/min stability for these d.fs with expressions analogous to (1.1) and (1.2) holding for all $j \in I_0$ are possible. The restriction $\alpha<1$, in the definitions of semi-Weibull, semi-Pareto and generalized semi-Pareto laws makes them mixtures of exponentials (now $R(x)$'s are LTs) and $x$ to $j \in I_0$, yield their discrete versions as in (2.2). Thus geometric-max or min stability characterizes discrete semi-Pareto and discrete versions of Propositions.2.1 and 2.2 characterizes geometric (on $I_0$), Sibuya and discrete semi-Weibull. That of Proposition.2.3 with $\alpha<1$, $ap = 1 = \beta k$ and $p = c^\alpha$ is:

**Proposition.2.7** A d.f on $I_0$ is Harris-min stable iff it is discrete generalized semi-Pareto.

## ACKNOWLEDGEMENT

Authors thank the referee for the suggestions.